\documentclass{gtart}
\usepackage{amssymb,amsmath,amsthm,epsfig,graphicx,color,psfrag}

\newtheorem{thm}{Theorem}%[section]
\newtheorem{prop}[thm]{Proposition}
\newtheorem{lemma}[thm]{Lemma}

\theoremstyle{definition}
\newtheorem{defn}[thm]{Definition}

\newtheorem*{rem}{Remark}
\newtheorem{remnum}[thm]{Remark}
\newtheorem{exam}[thm]{Example}

\newtheorem{keyobs}[thm]{Key Observation}

\def\<{\langle}\def\>{\rangle}
\def\wtil{\widetilde}

\def\co{\colon\thinspace}
\def\inv{^{-1}}

\def\graph{\Gamma} % notation for the defining graph of the RAAG
\def\AG{A} % notation for the RAAG

\def\square{{\vcenter{\hrule height.4pt
      \hbox{\vrule width.4pt height5pt \hskip5pt
           \vrule width.4pt}
      \hrule height.4pt}}}

\begin{document}
\title{\textbf{The conjugacy problem in right-angled Artin groups and 
their subgroups}}
\author{John Crisp, Eddy Godelle, Bert Wiest}
\address{Institut de Math\'emathiques de Bourgogne, CNRS UMR 5584 \\
Universit\'e de Bourgogne, B.P. 47870, 21078 Dijon cedex, France\\ 
{\tt jcrisp@gmail.com}
\\
and
\\
Laboratoire de Math\'ematiques Nicolas Oresme,
CNRS UMR 6139\\
Universit\'e de Caen, 14032 Caen cedex, France\\
{\tt eddy.godelle@math.unicaen.fr}
\\
and
\\
IRMAR, Campus de Beaulieu, CNRS UMR 6625\\ 
Universit\'e de Rennes 1, 35042 Rennes, France\\
{\tt bertold.wiest@univ-rennes1.fr}}\\
\vspace{6mm}

\begin{abstract}
We prove that the conjugacy problem in right-angled Artin groups (RAAGs),
as well as in a large and natural class of subgroups of RAAGs, can be 
solved in linear-time. This class of subgroups contains, for instance,
all graph braid groups (i.e.~fundamental groups of configuration 
spaces of points in graphs), many hyperbolic groups, and it 
coincides with the class of fundamental groups of ``special
cube complexes'' studied independently by Haglund and Wise.
\end{abstract}
\keywords{
right-angled Artin group, partially commutative group, graph group,
graph braid group, conjugacy problem, cubed complex, special
cube complex.}

\primaryclass{20F36}
\secondaryclass{20F10, 20F65.}
\makeshorttitle

%|<------------------------------------------------------------------------>|

\section{Introduction}

It is well known that the conjugacy problem in free groups can be solved
in linear-time by a RAM (random  access memory) machine. 
This result has been generalized in two different directions. On the one 
hand, Epstein and Holt~\cite{EpHo} have shown that the conjugacy problem 
is linear in all word-hyperbolic groups. 
On the other hand, Liu, Wrathall and Zeger have proved the analogue result
for all right-angled Artin groups (\cite{LWZ}, based on~\cite{Wr}).
Note that these groups are also called ``partially commutative groups''
or ``graph groups'' in the literature.

The aim of the present paper is to extend the second approach, in order to
prove linearity of the conjugacy problem in a large class of subgroups
of right-angled Artin groups. Very roughly speaking, the subgroups in 
question are fundamental groups of cubical complexes, sitting inside 
the right-angled Artin group in a convex fashion.
%, namely the isometrically embedded ones
This class of groups has previously been studied by Crisp and Wiest
\cite{CW1,CW2}, and independently by Haglund and Wise~\cite{HagWise},
as fundamental groups of so-called {\it special cube complexes} (or,
more precisely, $A$-special cube complexes). 

The class of groups considered in this paper contains in particular all 
graph braid groups \cite{Ab,AbGh,FarSab1,FarSab2,Sab} and more generally 
all state complex groups \cite{AbGh2,GhrPet}. These classes of groups 
have attracted 
considerable interest recently, which stems partially from their close 
relations to  robotics \cite{AbGh,GhrPet}. Indeed, our results can be 
interpreted as giving very efficient algorithms for motion planning
of periodic robot movements. However, our results also apply to
the various word-hyperbolic groups discussed in~\cite{CW1,CW2} -- in
particular, to all surface groups except the three simplest non-orientable
ones.

The present paper raises the stakes on the conjecture of Haglund and 
Wise~\cite{HagWise} that all Artin groups
(e.g. braid groups) are virtually fundamental groups of special cube 
complexes. If this conjecture was known to be true, then our work would 
imply that Artin groups have finite index subgroups where
the conjugacy problem can be solved in linear time. 

The plan of the paper is as follows.
In the second section we present an alternative approach
to the conjugacy problem in right-angled Artin groups, different from the
one of Liu, Wrathall and Zeger, but rather close in spirit to the methods
of Lalonde and Viennot~\cite{Lalonde,Vie}.
In the third section we prove that isometrically embedded 
subgroups of right-angled Artin group inherit a linear-time solution to the
conjugacy problem from their supergroups.

%|<------------------------------------------------------------------------>|

\section{The conjugacy problem in RAAGs is linear-time}
\label{S:RAAGAlg}

We recall that a right-angled Artin group is a group given by a
finite presentation, where every relation states that some
pair of generators commutes. Graphically, a right-angled Artin group
$\AG$ can be specified by a simple graph~$\graph_\AG$, where the generators
of~$\AG$ correspond to the vertices of~$\graph_\AG$, and a pair of
generators commutes if and only if the corresponding vertices are
\emph{not} connected by an edge. 
Note that the opposite convention (connecting \emph{commuting} generators 
by an edge) is also very common, but in the present paper we shall stick to 
this convention.

Right-angled Artin groups have been widely studied in the last 
decades -- see~\cite{Cha} for an excellent survey. 
Several solutions to the word and conjugacy problem have been found. 
It seems to be difficult to have a complete bibliography of the large 
number of articles on these two problems. The first solutions 
to the word and the conjugacy problem was obtained by 
Servatius in~\cite{Ser}. In~\cite{VWy}, Van Wyk constructed a normal form 
in right-angled Artin groups and proved that these groups are biautomatic. 
Indeed, even thought our point of view is very different from Van Wyk's,
the normal form constructed in the present paper is the very similar to his.
One of the main papers regarding the algorithmic complexity of
these two problems is~\cite{LWZ} (based on~\cite{Wr}) by Liu, Wrathall and 
Zeger, which proves that they are both of linear complexity.

The word problem in partially commutative \emph{monoids} has also been 
widely studied and numerous papers appeared on that topics. Several 
approaches appeared to be successful. In~\cite{CaF}, Cartier and Foata 
constructed a normal form on  partially commutative monoids, 
and then obtained the first solution to the word problem. This normal form 
is the restriction of the normal form obtained in~\cite{VWy}. More recently, 
Viennot introduced in~\cite{Vie} a new tool, the so-called \emph{Viennot's 
piling}, based on a geometrical representation of partially commutative 
monoids. Several works deal with this tool (see for instance~\cite{Dub1} 
and~\cite{Lalonde}). The Viennot piling method associates a piling to 
each element of a partially commutative monoid and thereby provides a 
linear-time solution to the words problem in such a monoid. As 
remarked by Krob, Mairesse, and  Michos in~\cite{KMM}, this piling is 
canonically related to the normal form constructed in~\cite{CaF}. In 
\cite{Lalonde}, Lalonde introduces and uses the notion of a \emph{pyramid} 
in order to study the conjugacy problem in partially commutative monoids. 
In the present paper, we are going to extend the notions of a piling and 
of a pyramid to the context of right-angled Artin \emph{groups}, and use 
them in order to obtain a linear-time solution (to the word problem and) to the conjugacy problem. This leads us to introduce the notion of a cyclic normal form.

In order to get an intuition for the nature of the conjugacy problem in
right-angled Artin groups, let us first consider the relatively easy
case of free groups. Given two cyclic words of length~$\ell_1$ 
and~$\ell_2$ respectively, there is a two step algorithm which can
be performed in time $O(\ell_1+\ell_2)$ on a RAM machine: first each word 
can be cyclically reduced in time~$O(\ell_1)$ and~$O(\ell_2)$, respectively.
If the reduced words have different lengths, then they are not conjugate.
If they have the same length~$\ell$, then they can be compared in time
$O(\ell)$ using standard pattern matching algorithms,
like the Knuth-Morris-Pratt algorithm, the Boyer-Moore algorithm, or 
algorithms based on suffix-tree methods 
-- see \cite{KMP,BM,AHU,Gus,Ste}. It should be stressed that on a 
Turing machine these algorithms take time~$O(\ell \log(\ell))$.

In the sequel, we assume that~$\AG$ is a fixed right-angled Artin group 
given by a fixed presentation. We denote by~$\{a_1,\cdots, a_N\}$ the 
generating set of~$A$ associated with this presentation.

The aim of this section is to provide an algorithm which does, very 
roughly speaking, the following: given a word~$w$, another word~$w'$ 
with smaller or equal length is created in linear time such 
that~$w$ and~$w'$ represent conjugate elements of~$\AG$. Furthermore, 
the word~$w'$ depends only on the conjugacy class in~$\AG$ of the 
element represented by~$w$, up to a cyclic permutation of its letters. 
This yields a linear-time solution to the conjugacy problem in~$\AG$ 
because, given words $w$ and $v$ we can compute the canonical
cyclic words $w'$ and $v'$ representing their conjugacy classes,
and compare those by one of the algorithms mentioned above.
%the algorithm does not increase the length and, as already mentioned, 
%there are standard linear time algorithms for deciding whether or not 
%two cyclic words are equal.

%............................................................................

\subsection{The word problem is linear-time} 
 We start by recalling the following classical lemma.
\begin{lemma}\cite{Ser}\label{L:reduced}
Any element of~$\AG$ can be represented by a \emph{reduced}
word (one which does not contain a subword of the form 
$a_i^{\pm 1} x a_i^{\mp 1}$, where all letters of~$x$ commute with~$a_i$).
Moreover, any two reduced representatives of the same element are related by a finite number
of commutation relations -- no insertions/deletions of trivial pairs
are needed.
\end{lemma}

Now we introduce our main tool, the notion of a \emph{piling}.

\begin{defn}
An \emph{abstract piling} is a collection of~$N$ words, one for each 
generator~$a_i$ of~$\AG$, over the alphabet with three symbols~$\{+,-,0\}$. 
\end{defn}
 The word associated with the generator~$a_i$ will be called the \emph{$a_i$-stack} of the abstract piling.
The product of two abstract pilings is defined as the piling obtained by concatenation of the corresponding stacks.\\

We define a function~$\pi^\star$ on the 
set~$\{a_1^{\pm 1},\ldots,a_N^{\pm 1} \}^*$ of words on the~$2n$ 
letters~$a_1,a_1\inv,\ldots,a_N,a_N\inv$ that associates an abstract 
piling to every word in the following way: starting with the empty piling,
we read the word from left to right. When a letter~$a_i^\epsilon$  
is read, we check what the last letter of the~$a_i$-stack of the piling is. 
If this letter is different from~$-\epsilon$ (the no-cancellation cases: 
the~$a_i$-stack is empty, or finishes either with~$0$ or~$\epsilon$), then 
we append a letter~$+$ or~$-$ at the end of the~$a_i$-stack of the piling (the sign of~$\epsilon$). Moreover, we also append a letter~$0$ at the end 
of each of the~$a_j$-stacks associated with a generator~$a_j$ that does 
not commute with the generator~$a_i$. On the other hand, if the last 
letter of the~$a_i$-stack is $-\epsilon$ (the cancellation case), 
then we erase this last letter, and we also erase the terminal letter 
of each of the~$a_j$-stacks of the piling associated with
a generator~$a_j$ that does not commute with the generator~$a_i$ --
note that the terminal letter of the~$a_j$-stack is necessarily ``$0$''. 

\begin{defn}\label{D:piling}
A \emph{piling} is an abstract piling in the image of the 
function~$\pi^\star$. The set of pilings is denoted~$\Pi$.
\end{defn}

We observe that the number of letters~$+$ and~$-$ occuring in the piling
$\pi^\star(w)$ is at most equal to the length of the word~$w$. Moreover, 
it is immediate from the description of the function~$\pi^\star$ that, 
given a word~$w$ of length~$\ell$, the piling~$\pi^\star(w)$ can be 
calculated in time~$O(\ell)$ (linear-time).

It may be helpful to keep in mind the following physical interpretation 
of a piling: we have~$N$ vertical sticks, labelled by the 
generators~$a_1,\ldots, a_n$, with beads on it; the beads are 
labelled by~$+$,~$-$ or~$0$ such that 
when reading from bottom to top the sequence of labels of the beads on 
the~$a_i$-stick, we obtain the~$a_i$-stack of the piling. A letter~$a_i$ 
or~$a_i\inv$  of the word~$w$ corresponds to a set of beads (which we call 
a \emph{tile}), consisting of one bead labelled~$+$ or~$-$ on the corresponding stick, and one bead labelled~$0$ on each of the 
sticks corresponding to generators of~$\AG$ which do not commute with $a_i$;  
each~$0$ labelled bead is connected to the~$\pm$~labelled bead by a thread. 
The rule is: on a stick, adjacent~$0$-beads  can commute with (``slide 
through'') each other, but~$0$-beads do not commute with~$\pm$-beads. 
In this physical model, we construct the image of a word by adding beads 
from the top, and removing opposed tiles when one obtains on a stick two 
adjacent~$\pm$-beads with opposite signs. In fact, when we are dealing
with the word problem we can forget about the threads between the beads,
but they are helpful for thinking about the conjugacy problem. 

\begin{exam} \label{E:exemplesuivi1} In the group~$A$ with group 
presentation
$$\langle a_1,a_2,a_3,a_4\mid a_1a_4 = a_4a_1\ ;\ a_2a_3 = a_3a_2\ ;\ a_2a_4 = a_4a_2 \rangle$$
we can calculate the piling~$p$ of the word
$a_2^{-2}a_4^{-1}a_3a_2a_4a_1a_2a_1^{-1}a_2^2a_4^{-1}$
as indicated in Figure~\ref{F:calculatepiling}. 
\begin{figure}[ht]
\centerline{\epsfig{file=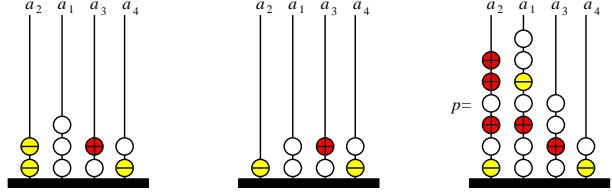,height= 2.5cm}}
\caption{The pilings of the prefixes~$a_2^{-2}a_4^{-1}a_3$ and
$a_2^{-2}a_4^{-1}a_3a_2$, and of the full word
$a_2^{-2}a_4^{-1}a_3a_2a_4a_1a_2a_1^{-1}a_2^2a_4^{-1}$}
\label{F:calculatepiling} 
\end{figure}
\end{exam}
The map~$\pi^\star$ induces a well-defined function~$\pi\co \AG\to \Pi$
 because words
representing the same element of~$\AG$ have the same image under~$\pi^\star$: the image of a word is unchanged 
by applying a commutation relation, and by inserting or deleting a trivial pair~$a_ia_i\inv$ or~$a_i\inv a_i$.  Now, from the definitions it is immediate that no cancellation occurs during the construction of the piling~$\pi^\star(w)$ of a reduced word~$w$. Then,  the identity of $A$ is the unique element of $A$ whose image by~$\pi$ is the trivial piling, and therefore the word problem is solved in linear-time: a word~$w$ represents the identity if and only if its piling~$\pi^\star(w)$ is trivial; this piling can be built in linear-time. 

The following notion will be extremely useful in the next section when we consider the conjugacy problem.

\begin{defn}
Let $w$ be a reduced word.\\
(i) We say that $w$ is \emph{initially normal} when $w$ is trivial or when the index of its first letter is greater or equal to the index of the first letter of any equivalent reduced word.\\
%\marg{replaced ``factor'' by ``suffix''.  B.}
(ii) We say that $w$ is \emph{normal} when all its suffixes are initially normal.
\end{defn}
We remark that all the factors of a normal word are normal words.
\begin{prop}\label{P:unnmfm}
 Any element of $A$ has a unique normal reduced representative word.
\end{prop}
\begin{proof} 
For any reduced word~$w = a^{\varepsilon_1}_{i_1}\cdots a^{\varepsilon_k}_{i_k}$, where $\varepsilon_j = \pm 1$, we set 
$$
\Omega(w) = \{(r,s)\mid 1\leq r<s\leq k\textrm{ and }i_r<i_s\}.
$$
Let~$a$ be in~$A$. In order to prove that $a$ has normal reduced
representative word, we choose, among all words representing $a$,
a word $w$ for which the number $\#\Omega(w)$ is as small as possible
(possibly equal to zero). This word $w$ is minimal.

We shall prove uniqueness of the normal representative by induction on
the length. If $a$ is of length $1$, i.e.\ if $a=a_i^\varepsilon$
for $\epsilon=\pm 1$, then uniqueness is obvious. 

Now suppose that $a$ 
has two normal reduced representatives 
$w=a^{\varepsilon_1}_{i_1}\cdots a^{\varepsilon_k}_{i_k}$ and
$w'=a^{\varepsilon'_1}_{i'_1}\cdots a^{\varepsilon'_k}_{i'_k}$.
Since the the suffixes of length $k-1$ of $w$ and $w'$ are again
normal, it is, by induction hypothesis, sufficient to prove that~$a_{i_1}^{\varepsilon_1} = a^{\varepsilon'_1}_{i'_1}$.   
Since $w$ and $w'$ are normal, we have $i_1 = i'_1$. Now, the exponents 
also have to be equal by Lemma~\ref{L:reduced}: we can not transform 
the word~$a^{\varepsilon_1}_{i_1}\cdots a^{\varepsilon_k}_{i_k}$ into 
the word~$a^{-\varepsilon_1}_{i_1}a^{\varepsilon'_2}_{i'_2}\cdots 
a^{\varepsilon'_k}_{i'_k}$ by using commutation relations only: 
starting from the reduced~$w$, no word of the form 
$ua^{\varepsilon_1}_{i_1}a^{-\varepsilon_1}_{i_1}u'$ can appear by 
any sequence of commutation relations.
\end{proof}

In the sequel, we call this unique normal reduced word representing~$a$ 
the \emph{normal form} of~$a$.

\begin{prop}~\label{P:existsigmamodif}
There is a linear-time algorithm that associates to each piling~$p$ a normal
word~$\sigma^\star(p)$ such that~$\pi^\star(\sigma^\star(p)) = p$. Furthermore, for any element~$a$ of~$A$ the word~$\sigma^\star(\pi(a))$ is the normal form of $a$.
\end{prop}

\begin{exam}\label{E:exemplesuivi2} Using the notation of 
Example~\ref{E:exemplesuivi1}, the word~$\sigma^\star(p)$ is equal to~$a_4^{-1}a_3$ $a^{-1}_2a_1a_2a_1^{-1}a_2a_2$. The calculation is shown in 
Figure~\ref{F:extract}.
\begin{figure}[ht]
\centerline{\epsfig{file=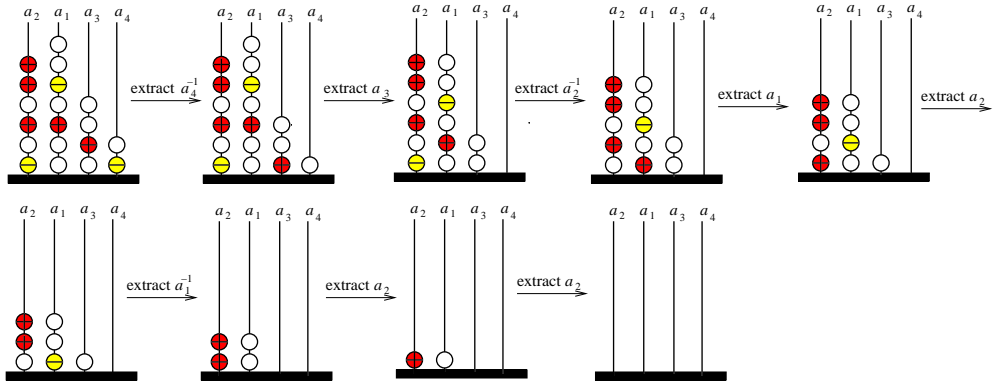,height= 5cm}}
\caption{The word $\sigma^\star(p)$ associated to a piling $p$}
\label{F:extract}
\end{figure}
\end{exam}
\begin{proof}[Proof of Proposition~\ref{P:existsigmamodif}]
Let~$p$ be a piling. By definition, this means that there exists an
element $a$ of~$A$ such that $\pi(a)=p$. In order to prove 
Proposition~\ref{P:existsigmamodif}, it suffices to find an algorithm
for constructing in linear time a word~$\sigma^\star(p)$, and to prove
that $\sigma^\star(p)$ is a normal reduced representative of $a$.

We start with the observation that the element $a$ has a reduced
representative starting with the letter $a_i^{\pm 1}$ if and only if
the $a_i$-stack of the piling is nonempty and starts with the letter
$+$ or $-$, respectively (not with the letter $0$). 

We associate to~$p$ a normal reduced word~$\sigma^\star(p)$ by induction 
on the number of letters $+$ and $-$ in~$p$ in the following way. 
If~$p$ is empty then~$\sigma^\star(p)$ is the empty word. Otherwise, 
let~$i$ be the largest index with the property that the~$a_i$-stack 
of~$p$ is nonempty and starts with the letter~$+$ or~$-$, not with~$0$. 
Then, according to this sign, we define the first letter 
of~$\sigma^\star(p)$ to be~$a_i$ or~$a_i^{-1}$, respectively. Then we 
remove the tile consisting of the first letter ($+$ or~$-$) 
of the~$a_i$-stack, and of the initial letter (which has to be $0$)
of each of the~$a_j$-stacks associated with a generator~$a_j$ that does 
not commute with~$a_i$.  
What remains is a piling~$p_1$ with strictly fewer letters. Thus the
word $\sigma^\star(p_1)$ is already defined, by induction hypothesis,
and we define the word $\sigma^\star(p)$ by concatenation 
$\sigma^\star(p)=a_i^{\pm 1} \sigma^\star(p_1)$.

We claim that the word $\sigma^\star(p)$ is a normal reduced 
representative of~$a$; indeed, in the above construction we see that
the first letter of $\sigma^\star(p)$ is also the first letter of
some reduced representative of $a$. By induction, the whole word
$\sigma^\star(p)$ is a reduced representative of $a$. Moreover, 
the word $\sigma^\star(p)$ is initially normal, by construction, 
and by induction its suffix $\sigma^\star(p_1)$ is normal. Hence the
whole word $\sigma^\star(p)$ is normal.
\end{proof}

\subsection{Cyclic normal forms and pyramidal pilings}\label{SS:CyclicNF}
We are now ready to attack the conjugacy problem.  

\subsubsection{Cyclically reduced words and cyclically reduced pilings}

We recall that a \emph{cycling} of a reduced word~$w$ is the operation of
removing the first letter of the word, and placing it at the end
of the word. A word is called \emph{cyclically reduced}
if it is reduced and if any 
word obtained from it by a sequence of cyclings and commutations is 
still reduced  -- in other words, if it is not 
of the form~$x_1 a_i^{\pm 1} x_2 a_i^{\mp 1} x_3$, where all the letters of
$x_1$ and~$x_3$ commute with~$a_i$. As far as we know, all known solutions to the conjugacy problem in RAAGs are based on the
following lemma.

\begin{lemma}
Two cyclically reduced words represent conjugate elements of~$\AG$ if
and only if they are related by a sequence of cyclings and commutation
relations.
\end{lemma}

Therefore two reduced words~$w_1,w_2$ with letters in 
$\{a_1^{\pm 1},\ldots,a_n^{\pm 1}\}$ represent conjugate elements of~$\AG$ 
if and only if there is a sequence of words
$$
w_1\stackrel{\mathrm{red}}{\longrightarrow} v_1 \leftrightarrow v_2 \stackrel{\mathrm{red}}{\longleftarrow} w_2 
$$
where the two arrows labelled ``red'' represent two sequences of cyclic 
reductions down to cyclically reduced words and  the 
arrow~$\leftrightarrow$ represents a finite sequence of cyclings and commutation relations.

\begin{defn}
If, in a piling~$p$, the~$a_i$-stack starts ({\it resp.} finishes) 
with a letter~$+$ or~$-$, the \emph{bottom} \emph{$a_i$-tile} ({\it resp.} the \emph{top} \emph{$a_i$-tile}) of~$p$ is the \emph{sub-piling} formed by the first ({\it resp.} last) letter of the~$a_i$-stack and the first ({\it resp.} last) letter of the~$a_j$-stacks such that~$a_i$ and~$a_j$ do not commute in~$\AG$.
\end{defn}

\begin{exam} \label{E:exemplesuivi3} With the notation of Example~\ref{E:exemplesuivi1}, Figure~\ref{F:topbot} gives an example of top 
and bottom tiles of a piling.
\begin{figure}[ht]
\centerline{\epsfig{file=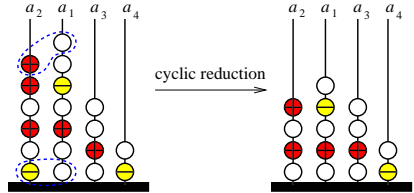,height= 2.5cm}}
\caption{a top~$a_2$-tile and a bottom~$a_2$-tile, and the associated cyclic reduction} \label{F:topbot}
\end{figure}
\end{exam}

\begin{defn}
If in a piling~$p$ the~$a_i$-stack starts with the letter~$+$ and ends 
with~$-$, or vice versa, a \emph{cyclic reduction} is the act of 
removing both top and bottom~$a_i$-tiles. We say that the piling is 
\emph{cyclically reduced} if no cyclic reduction is possible. 
\end{defn}

Note that cyclically reducing a piling yields again a piling.
We remark that there is an obvious linear-time algorithm for transforming 
any piling into a cyclically reduced one by a finite sequence of cyclic 
reductions. We also observe that for a reduced 
word~$w\in \{a_1^{\pm 1},\ldots,a_N^{\pm 1}\}^*$, cycling of~$w$ corresponds 
to a cycling of its piling, and that~$w$ is cyclically reduced if and only 
if the piling~$\pi^\star(w)$ is.

Now we have a fast algorithm for cyclically reducing words and pilings.
In contrast to the case of free groups, however, the reduced words which
we can obtain are \emph{not} unique up to cyclic permutation.
In order to circumvent this problem, we shall introduce in the sequel the notion of a \emph{cyclic normal form}.

\subsubsection{Non-split words and non-split pilings}\label{SSS:nonsplit}

Our first objective is to restrict the conjugacy problem to the case 
of \emph{non-split cyclically reduced words} ( or \emph{pilings}). 
We recall that a graph~$\Gamma_A$ is associated to the right-angled 
Artin group~$A$. 

\begin{defn}
Let~$w$ be a reduced word different from~$1$, and let~$p$ 
be its image by~$\pi^\star$. Consider~$\Delta(p)$ (or~$\Delta(w)$) the full 
subgraph of~$\Gamma_A$  whose vertices are those whose correponding stacks 
contain at least one bead different from~$0$ (in other words, the 
letters~$a_i$ such that~$a_i^{\pm 1}$ occurs in~$w$). Then, the word~$w$ 
and the piling~$p$ are said to be \emph{non-split} when the 
graph~$\Delta(p)$ is connected. 
\end{defn}

\begin{figure}[ht]
\centerline{\epsfig{file=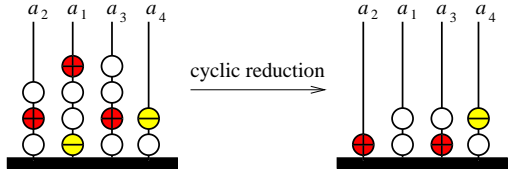,height=2.2cm}}
\caption{The word~$a_1^{-1} a_2 a_3 a_1 a_4^{-1}$ is not split, but
cyclic reduction yields a word which is split: 
$a_2(a_3a_4^{-1}) = (a_3a_4^{-1})a_2$} \label{F:splitpiling}
\end{figure}

In other words,~$w$ is non-split if and only if its set of letters cannot 
be separated in two disjoint subsets such that every letter of one of the 
subset commutes in~$A$ with every letter of the other subset. Clearly, 
it takes linear-time to obtain the set of vertices of the 
graph~$\Delta(p)$, and constant time (which depends on the 
graph~$\Gamma_A$) to decide if~$\Delta(p)$ is connected. If it is not, 
it takes still constant time to determine the connected 
components~$\Delta_1(w),\ldots, \Delta_k(w)$ of $\Delta(w)$. 
Figure~\ref{F:splitpiling} (which still uses the notation of 
Example~\ref{E:exemplesuivi1}) contains examples
of both split and non-split pilings. 

Now, if~$w$ is a cyclically reduced word that is split, then it 
is equivalent to a product~$w_1\cdots w_k$  of non-split cyclically 
reduced words, one for each connected component~$\Delta_i(w)$ of the 
graph~$\Delta(w)$; the graph~$\Delta_i(w)$ is equal to~$\Delta(w_i)$. 
Furthermore, once that the connected 
components~$\Delta_1(w),\ldots, \Delta_k(w)$ of $\Delta(w)$ are computed, 
appropriate words $w_1,\ldots, w_k$ can be obtained in linear-time.   

\begin{remnum}\label{R:ReductionToNonsplit}
The following observation will be crucial: if~$v$ is another cyclically 
reduced word, then then~$w$ and~$v$ represent conjugate elements if and 
only if two conditions are satisfied: 
firstly the graph~$\Delta(v)$ is equal to~$\Delta(w)$;
secondly, if $v_1,\ldots,v_k$ are words such that $\Delta(v_i) = \Delta_i(v)$ 
and such that $v$ is equivalent to the product $v_1\cdots v_k$, then
for each index~$i$ the words~$w_i$ and~$v_i$ represent conjugate elements. 
\end{remnum}

Therefore, in order to obtain a solution to the conjugacy problem in 
linear-time it is enough to consider the case of cyclically reduced
non-split words.

\subsubsection{Pyramidal piling and cyclic normal form}

To solve the conjugacy problem, we associate in the sequel  
a \emph{cyclic normal word} to each cyclically reduced non-split word. 
We first do the analogue of this in the framework of pilings: to each 
non-split cyclically reduced piling, we associate a \emph{pyramidal} 
piling. 

\begin{defn} Let~$p$ be a non-empty piling, and denote by~$i$ the 
smallest index such that the~$a_i$-stack contains an~$a_i^{\pm}$-bead.
We say that the piling~$p$ is \emph{pyramidal} if the first bead of 
every~$a_j$-stack except the~$a_i$-stack is either empty or
starts with the letter~$0$. In that case, we say that~$a_i$ is the 
\emph{apex} of the pyramidal piling.
\end{defn}
Note that a pyramidal piling has to be non-split.

\begin{lemma}\label{lemdecomp}
(i) Let~$p$ is a non-empty piling and denote by~$i$ the smallest index such that the~$a_i$-stack of~$p$ contains an~$a_i^{\pm}$-bead; then there exists a unique decomposition~$p_0\cdot p_1$ of~$p$ such that~$p_1$ is a pyramidal piling with~$a_i$ as apex, and~$p_0$ is a piling without~$a_i^{\pm}$-beads. 
Furthermore, one has 
the equality of words~$\sigma^\star(p) = \sigma^\star(p_0)\sigma^\star(p_1)$.

(ii) The above decomposition~$p_0\cdot p_1$ can be computed in linear-time on the number of~$\pm$-beads of the piling~$p$. 
\end{lemma}

\begin{exam} \label{E:exemplesuivi4} Using the notation of Example~\ref{E:exemplesuivi1}, Figure~\ref{F:decomp} gives an example of a decomposition of a piling.
\begin{figure}[ht]
\centerline{\epsfig{file=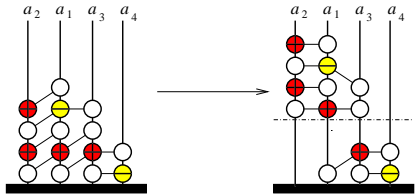,height=2.5cm}}
\caption{Decomposition of a piling as~$p_0\cdot p_1$}\label{F:decomp} 
\end{figure}
\end{exam}

\begin{proof}[Proof of Lemma~\ref{lemdecomp}]
We start by exhibiting a linear-time algorithm for finding such a 
decomposition of a given non-empty piling~$p$. 
Let~$p_0$ be the empty piling. Reading all the stacks (in the index order), 
obtain in linear-time the smallest index~$i$ for which the~$a_i$-stack 
contains a bead distinct from~$0$. Then, apply iteratively the following 
recipe: consider the largest index~$j$ (necessarily greater than~$i$) 
for which 
the~$a_j$-stack  starts with a letter~$+$ or~$-$; then remove all the 
beads in the bottom~$a_j$-tile, and add them to the top of the 
piling~$p_0$. When no more beads can be extracted from the bottom of the 
piling~$p$, then the construction of the factor~$p_0$ is complete, and what 
remains is the piling~$p_1$. This proves the existence part of (i), as well 
as part~$(ii)$. 
The formula~$\sigma^\star(p) = \sigma^\star(p_0)\sigma^\star(p_1)$ is 
now immediate by construction. For the uniqueness part of (i), we notice 
that in any decomposition~$p=p_0\cdot p_1$, the factor~$p_0$ has to 
contain exactly those tiles that can be extracted on the bottom from~$p$ 
without extracting any apex bead. 
\end{proof}

We call the piling~$p_0$ the~\emph{$0$-factor} of~$p$. 
Thus the piling~$p$ is pyramidal if and only if its~$0$-factor is empty.

In our physical interpretation, if~$i$ is the smallest index such that the~$a_i$-stack contains a~$a_i^\pm$-bead, we can lift up the first
$a_i^{\pm}$-bead along its stick to the first floor. Then some part of 
the piling stays on the ground, while some beads are lifted up. Here it 
is essential to keep in mind that each~$0$-bead is connected by a thread 
to a~$\pm$-bead, and that adjacent~$0$-beads on a stick can slide 
through each other. The factor that stays down is~$p_0$, the factor 
that is lifted up is~$p_1$. This latter factor has the structure of an 
upside-down pyramid supported by one of the apex-beads, hence the names.

\begin{figure}[ht]
\centerline{\epsfig{file=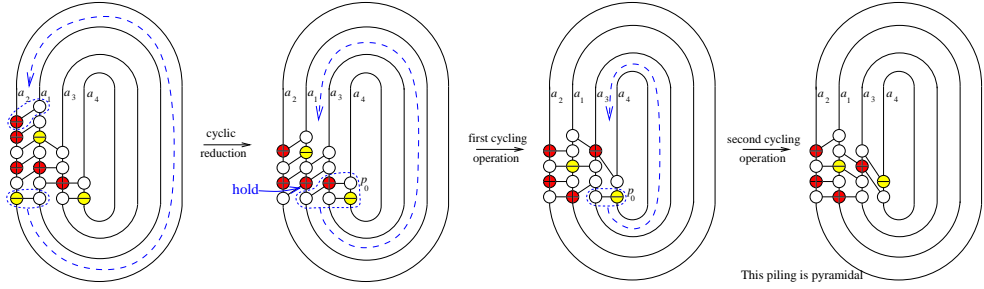,width=13cm}}
\caption{The calculation of a pyramidal piling} \label{F:CalcPyr}
\end{figure}

If in a cyclically reduced piling, the~$a_i$-stack 
starts with a letter~$+$ or~$-$, then one can perform a \emph{cycling}
of the bottom tile containing that bead to the top of the piling, i.e., 
one can move the initial letter~$+$ on the~$a_i$-stack, and the initial 
letters~$0$ on the stacks corresponding to letters that do not commute 
with~$a_i$, to the end of their respective stacks.
A physical interpretation (see Figure~\ref{F:CalcPyr}) of this procedure
is obtained by replacing the sticks by concentric hula hoops. A cycling 
of a bottom tile corresponds to the operation of cycling  the 
corresponding tile along the hula hoops.
 
\begin{prop}\label{P:KeyAlg}
There is an algorithm which takes as its input any non-split cyclically reduced piling~$p$ and which outputs a pyramidal piling that is obtained from the input piling by a finite sequence of cyclings. If the piling has~$\ell$ beads, then the algorithm requires~$O(\ell)$ cyclings, so its computational complexity is~$O(\ell)$.
\end{prop}
 
\begin{proof} 
The basic procedure of the algorithm is in two steps; given a cyclically reduced piling~$p$, first determine the 
$0$-factor~$p_0$ of the canonical decomposition (by the method of Lemma~\ref{lemdecomp}). 
Secondly, cycle all the tiles belonging to~$p_0$ in order to obtain a new piling.
This procedure takes time~$O(\ell)$.
The algorithm is simply to iterate this basic procedure until the 
factor~$p_0$ is empty. It remains to prove that there is a bound on the number of iterations which
depends only on the group~$\AG$, not on the piling~$p$. In fact,
if we denote~$i$ the smallest index such that~$p$ contains an~$a_i$-tile, and ~$\Delta(p)$ the full subgraph of the defining graph~$\Gamma$ defined
above, we claim that~$\max_{a_j\in\Delta(p)}\ dist_\Delta(a_i,a_j)$
is an upper bound on the number of iterations, where each edge of~$\Delta(p)$ has length~$1$. This quantity is
finite, because~$\Delta(p)$ is connected, and is bounded above by~$N$, the number of generators of the group~$\AG$ (which does not depend on the piling~$p$). This fact is obvious from the geometrical representation, and the proof is a straightforward induction: after the first iteration of the basic procedure, no~$a_j^{\pm}$-beads such that~$a_j$ is at distance~$1$ from~$a_i$ in~$\Delta(p)$  appear in the~$0$-factor; after a second iteration no~$a_j^{\pm}$-beads such that~$a_j$ is at distance at most~$2$ from~$a_i$ in~$\Delta(p)$  appear in the~$0$-factor, and so on. 
\end{proof}

Now, if~$w$ is a non-split reduced word, we can apply the above algorithm to the piling~$\pi(w)$ to obtain a pyramidal piling~$p$. Then, the words~$\sigma^\star(p)$ and~$w$ represent conjugate elements.

\begin{defn} Let~$w$ be a word in~$\{a_1^{\pm 1},\ldots,a_n^{\pm 1}\}^*$ that is reduced and cyclically reduced. We say that the word~$w$ is a \emph{cyclic normal form} if it is normal and all its cyclically conjugate words  are normal.
\end{defn}

Intuitively, if we regard~$w$ as a cyclic word, and we start reading
anywhere in the word, then the first letter that we read must always be
the largest-index letter that can be extracted on the left.
For instance, with the notation of Example~\ref{E:exemplesuivi1}, the word~$a_4^{-1}a_3a^{-1}_2a_1a_2a_1^{-1}a_2a_2$ is not a cyclic normal form: starting from the last letter and reading cyclically, we read out~$a_2a_4^{-1}\ldots$, which is already illegal, because the letters commute, and
$a_4^{-1}$ has a larger index than~$a_2$, so~$a_4^{-1}$ should come first. Another example: the word~$a_1a_2a_1^{-1}a_3a_4^{-1}a_2$ is a cyclic normal form. Our linear-time solution to the conjugacy problem is based on the two following results.

\begin{prop}\label{T:KeyAlg2} If $p$ is a non-split cyclically reduced 
pyramidal piling, then $\sigma^\star(p)$ is a cyclic normal form.
\end{prop} 

\begin{prop}\label{T:KeyAlg} Two cyclic normal forms represent conjugate elements if and only if they are equal up to a cyclic permutation. 
\end{prop}
\begin{proof}[Proof of Proposition~\ref{T:KeyAlg2}]
Firstly, we remark that a consequence of Lemma~\ref{L:reduced} is the following fact: if $a^\epsilon,b^\eta$ are letters ($\epsilon,\eta = \pm 1$) and $w$ is a reduced word such that~$b^{-\eta}w$ and~$wa^{\epsilon}$ are both reduced ({\it i.e.} no word equivalent to~$w$ starts and finishes with~$b^{\eta}$ and~$a^{-\epsilon}$, respectively) but the word~$b^{-\eta}wa^\epsilon$ is not reduced ({\it i.e.}~$wa^\epsilon$ is equivalent to some word that starts with~$b^\eta$), then~$a^\epsilon = b^\eta$ and all the letters of~$w$ commute with~$a$. Now, we know that~$\sigma^\star(p)$ is a normal cyclically reduced word. For a cyclically reduced word~$w$, the word~$ww$ is cyclically reduced (this follows directly from the above fact, or from the piling representation), and all the words cyclically conjugate to the former are subwords of the latter. Therefore, in order to prove the result, it is enough to prove that the word~$\sigma^\star(p)\sigma^\star(p)$ is normal. Assume that this is not the case. Since $\sigma^\star(p)$ is normal, we can then write $\sigma^\star(p) = w_1a_j^\eta w_2 = v_1a_i^\epsilon v_2$ such that~$a_j^\eta w_2v_1$ is initially normal but~$a_j^\eta w_2v_1a^\epsilon$ is not. In particular, there exists $a_k^\nu$, with $k>j$, such that $a_k^{-\nu} a_j^\eta w_2v_1$ is reduced but $a_k^{-\nu} a_j^\eta w_2v_1a_i^\epsilon$ is not. Since~$a_j^\eta w_2v_1a_i^\epsilon$ is  a subword of~$\sigma^\star(p)\sigma^\star(p)$, it is reduced. Using the above fact, we get that~$a_i^\epsilon = a_k^\nu$, and $a_k^\nu$ commute with all the letter of~$a_j^\eta w_2v_1$. In particular, the word~$\sigma^\star(p)$ is equivalent to~$a_k^\nu v_1v_2$. This is impossible because~$k$ is greater that $j$, and~$p$ is pyramidal. Therefore,~$\sigma^\star(p)\sigma^\star(p)$ is normal. 
\end{proof}
 
\begin{proof}[Proof of Proposition~\ref{T:KeyAlg}]  
The ``if'' implication is obvious, we have to prove the ``only if'' part.

Let~$w$ and~$w'$ be two cyclic normal forms that represent 
conjugate elements. Let~$i$ be the smallest index that appears in~$w$ and 
choose a distinguished letter~$a_i^\varepsilon$ in~$w$. As the words~$w$ 
and~$w'$ are cyclically reduced, there exists a sequence of 
words~$w_0 = w\to w_1\to\cdots w_r = w'$ that transforms~$w$ into~$w'$, 
such that~$w_{i+1}$ is obtained from~$w_i$ by a commutation 
or a cycling transformation. 

We can keep track of the distinguished 
letter~$a_i^\varepsilon$ along the transformations: 
write~$w_j = w'_ja_i^{\varepsilon}w''_j$. 
Assume the number~$\ell$ of commutations 
that involve the distinguished letter is positive. Since $w$ is a normal word, the first commutation~$w_j\to w_{j+1}$ that 
involves~$a^\varepsilon_i$ is ``from left to right'', i.e.\ it is of
the following 
form:~$w_j = w'_{j+1}a_{i'}^{\varepsilon'} a_i^\varepsilon w''_j$ 
and~$w_{j+1} = w'_{j+1}a_i^\varepsilon a_{i'}^{\varepsilon'} w''_j$ 
with~$i'>i$. 

Now, consider the last operation~$w_p \to w_{p+1}$ such 
that a letter~$a^\eta_k$ is exchanged with the distinguished 
letter~$a_i^\varepsilon$ from left to right: 
we have~$w_p = w'_{p+1}a_{k}^{\eta} a_i^\varepsilon w''_p$ 
and~$w_{p+1} = w'_{p+1}a_i^\varepsilon a_{k}^{\eta} w''_p$. 
We can also keep track of the distinguished letter~$a_{k}^{\eta}$. 
As long as the two letters do not cross each other again 
in the opposite direction, we have~$w''_qw'_q = y_q a_k^\eta z_q$ 
such that all the letters of $y_q$ commute with $a_k$ (where $q$ satisfies
$q>p$). In particular, 
$a_i^\varepsilon w''_qw'_q$ is not initially normal. But $w'$ is normal, 
so the two distinguished letters have to cross each other
again in the opposite direction: there exists $s$, with $p<s<r$, such 
that~$w_{s} = w'_{s}a_i^\varepsilon a_{k}^{\eta} w''_{s+1}$ 
and~$w_{s+1} = w'_{s}a_{k}^{\eta} a_i^\varepsilon w''_{s+1}$. Hence, 
we have a sequence
$$
w_p = w'_{p+1}a_{k}^{\eta} a_i^\varepsilon w''_p\to v'_{p+1}a_{k}^{\eta} a_i^\varepsilon v''_{p+1}\to\cdots\to v'_{s-1}a_{k}^{\eta} a_i^\varepsilon v''_{s-1} \to v'_{s}a_{k}^{\eta} a_i^\varepsilon v''_{s}\to w_{s+1}
$$ 
such that each word $v''_qv'_q$ is equal to the word $y_qz_qw'_q$. 
Thus we obtain a new sequence from $w$ to $w'$ with only $\ell-2$ 
commutations that involve the distinguished letter~$a_i^\varepsilon$. 

It follows that we can assume that no commutation involves the 
distinguished letter~$a_i^\varepsilon$ 
along the sequence~$w_0 = w\to w_1\to\cdots w_r = w'$. But this implies 
that the words~$a_{i}^{\varepsilon}w''_1w'_1$ 
and~$a_{i}^{\varepsilon}w''_rw'_r$ are equivalent. As they are both 
cyclic normal forms, they are normal words. 
Therefore they are equal by Proposition~\ref{P:unnmfm}. Hence, the 
words~$w$ and~$w'$ are equal up to a cyclic permutation.     
\end{proof}

Summing up, in order to decide whether two nonsplit cyclically 
reduced words represent conjugate elements, it suffices to decide 
whether their cyclic normal forms are equal (as cyclic words), and 
these cyclic normal forms can be calculated in linear time. More formally,
we have

\begin{thm} The conjugacy problem in a right-angled Artin group~$\AG$
is linear-time on the sum of the lengths of the two input words.
\end{thm}

\begin{proof} Here is a summary of the algorithm:\\
Given any two words $w$ and $v$,\\
(i)  produce the piling $\pi^\star(w)$, 
and then by cyclic reduction a cyclically reduced piling~$p$; similarly 
for the word~$v$ produce first the piling $\pi^\star(v)$, and cyclically 
reduce it to a piling~$q$;\\ 
(ii) factorize each of the pilings~$p$ and $q$ 
into non-split factors. If the collection of subgraphs $\Delta_i(p)$ 
and $\Delta_i(q)$ of the defining graph $\Gamma_\AG$ do not coincide,
output ``NO, $w$ and $v$ do not represent conjugate elements'' and stop. 
Otherwise,\\
(iii) if $p=p^{(1)}\cdot\ldots\cdot p^{(k)}$ and
$q=q^{(1)}\cdot\ldots\cdot q^{(k)}$ are the factorizations found in
step (ii), then for $i=1,\ldots,k$ do the following\\
\begin{enumerate}
\item[(a)]  transform the non-split cyclically reduced pilings~$p^{(i)}$ 
and~$q^{(i)}$ into pyramidal pilings $\wtil p^{(i)}$ and $\wtil q^{(i)}$, 
using a sequence of cyclings. Then produce the words in cyclic normal 
form $\sigma^\star(p^{(i)})$ and $\sigma^\star(q^{(i)})$;
\item[(b)] decide whether the words in cylic normal form found in the 
previous steps are the same up to cyclic permutation 
(in linear-time, using a standard algorithm). If they are not, 
answer ``NO'' and stop.
\end{enumerate}
(iv) answer "YES".
\end{proof}

\subsection{Calculating the centralizer of an element}\label{SS:centralizer}

The centralizer of a cyclically reduced element of $\AG$ has a canonical 
finite generating set:
suppose that $w$ is a cyclically reduced word, written as a product
of cyclically reduced non-split words $w=w_1\cdots w_k$, c.f.\ 
Section~\ref{SSS:nonsplit}. Then, according to~\cite{Baudisch2}, for each 
$i$ in $\{1,\ldots,k\}$ there exists a unique maximal infinite-cyclic 
subgroup of $\AG$ containing $[w_i]$, generated by some cyclically
reduced element~$[z_i]$, 
and by~\cite{Ser} the centralizer of $[w]$ in~$\AG$ is generated by
\begin{enumerate}
\item the elements $[z_i]$, and
\item the generators of $\AG$  which commute with all the generators 
occurring in~$w$. 
\end{enumerate}

In the next section we will need to algorithmically determine 
explicit representatives of these generators, in the special
case where the words $w_i$ are cyclic normal forms. 

\begin{prop}\label{P:commutator}
There is a linear-time algorithm which takes as its input
a cyclically reduced word $w$, decomposed as a product
of words in cyclic normal form $w=w_1\cdots w_k$, and which outputs
the canonical generating set of the centralizer of $w$.
\end{prop}

\begin{proof}
It takes linear time to determine the graph $\Delta(w)$, and then 
constant time to deduce from this the generators of type (2).

Now we turn to the generators of type (1), i.e.\ the minimal roots
$[z_i]$ of the elements $[w_i]$.
As a first step, we claim that periodicity of elements is visible in 
their cyclic normal form. More precisely, if one of the words~$w_i$
is equivalent to a word of the form~$\wtil{z}_i^{\,r}$ for some 
word~$\wtil{z}_i$ and some integer~$r$, then the word~$w_i$ itself is 
of the form~$z_i^r$, for some word $z_i$. 
In order to prove this claim, we observe that 
$\wtil{z}_i$ is equivalent to a word $z_i$ in cyclic normal form 
(because the $0$-factor of $p(\wtil{z}_i)$
must divide the $0$-factor of $p(w_i)$, which is the trivial word).
Now the word $z_i^r$ is still in cyclic normal form (c.f.\ the proof
of Proposition~\ref{T:KeyAlg2}), and it is equivalent to $w_i$. 
Therefore we have $z_i^r=w_i$.

We claim that for each of the factors $w_i$, the desired minimal 
root~$z_i$ of~$w_i$ is detectable in linear-time: we can calculate a 
pair $(z_i,r)$, where $z_i$
is a word and~$r$ an integer with $z_i^r=w_i$, and~$r$ is maximal
among all such pairs. Indeed, this algorithm works as follows:
consider the word $w_i^*$ obtained by removing the 
first letter from the word $w_i w_i$. Then find the starting point of the 
first occurrence of $w_i$ as a subword of $w_i^*$ -- this can be done by 
standard algorithms, like the Boyer-Moore algorithm, in time $O(\ell_i)$,
where $\ell_i$ denotes the length of $w_i$. 
If this starting point is at the $\ell_i$th letter of $w_*$, then there 
is no periodicity. If on the other hand the starting point is at the 
$t$th letter with $t<\ell_i$, then let $z_i$ be the prefix of~$w_i$ 
of length~$t$. By construction we have an equality of words 
$z_i w_i=w_i z_i$. This implies that the words $w_i$ and $z_i$ have a 
common root. By the choice of $z_i$, this root has to be $z_i$ itself
and for $r:=\ell_i/t$ we have an equality of words $w_i=z_i^r$.
Finally, by the choice of $t$, no prefix of $w_i$ of length less than~$t$
can be a root of~$w_i$, so~$z_i$ is indeed the minimal root.
\end{proof}

%|<------------------------------------------------------------------------>|

\section{The conjugacy problem in subgroups of RAAGs}

In the previous section we saw that the conjugacy problem in a fixed 
right-angled Artin group can be solved in linear-time on a RAM-machine 
with constant that depends only on the group. 
In this section we shall prove analogue results for a large class of 
subgroups of right-angled Artin groups, namely those considered in the 
papers~\cite{CW1,CW2}, as well as in~\cite{HagWise}. 

\subsection{A class of subgroups of RAAGs}

Every right-angled Artin group~$\AG$ admits a 
finite~$K(\AG,1)$, called the Salvetti complex of $\AG$, which
we shall denote~$Y$ and which can be constructed explicitely from the
presentation of~$\AG$. It is a cubed complex which has one 
single vertex, and one edge of length~$1$ for every generator of~$\AG$. 
Moreover, for every~$n$-tuple of mutually commuting generators of~$\AG$, 
there is one $(n+1)$-torus in~$Y$. 
We equip every cell, of any dimension, of this complex with the flat 
metric, in the sense that in the universal cover~$\wtil{Y}$ every cell 
is a Euclidean cube of sidelength~$1$. Then the complex is locally CAT(0), 
and its universal cover~$\wtil{Y}$ is CAT(0). For instance, for the group 
$A=\mathbb Z^2=\langle a_1,a_2 \ | \ [a_1,a_2]=1\rangle$,
the complex~$Y$ is a torus, constructed out of one vertex, two edges,
and one square which glued to the 1-skeleton according to the commutation
relation. See~\cite{CW1} for details. The reader should note that as soon 
as an orientation is chosen on each edge ({\it i.e.} simple loop) of~$Y$, 
one obtains an explicit isomorphism between~$A$ and~$\pi_1(Y)$ such that 
the image of each generator~$a_i$ of~$A$ is represented by the simple loop 
labelled by~$a_i$ traversed in the positive direction.

Now, suppose that~$X$ is a finite locally $CAT(0)$ cubed complex, and 
consider a cubical map~$\Phi\co X\to Y$, sending each open cube of 
$X$ bijectively and locally isometrically to a cell of the same dimension 
in~$Y$. (Here $Y$ still denotes the Salvetti complex of some 
right-angled Artin group.) If one of the vertices of $X$ is designated
as its basepoint, then such a mapping induces 
a homomorphism $\Phi_*\co \pi_1(X) \to \pi_1(Y)$. See 
Figure~\ref{F:counterex} for an example where~$X$ and~$Y$ are 
1-dimensional complexes. 

We need some more notation: for any
vertex $x$ of $X$, we denote by $\Phi_{lk}\co lk(x,X)\to lk(\Phi(x),Y)$ 
the induced map from the link of~$x$ in~$X$ to the link of~$\Phi(x)$ in~$Y$
We shall be interested in the following two properties which our map $\Phi$
may have:
\begin{itemize}
\item The convexity property: for any vertex~$x$ of~$X$, and any two 
  vertices of~$lk(\Phi(x),Y)$ which belong to the image
 ~$\Phi_{lk}(lk(x,X))$ and which are connected by an edge, the 
  connecting edge belongs to the image 
 ~$\Phi_{lk}(lk(x,X))$, as well. 
\item The injectivity property: the map of universal covers
 ~$\wtil{\Phi}\co \wtil X \to \wtil Y$ is injective. In particular,
 ~$\Phi_*\co \pi_1(X) \to \pi_1(Y)$ is a monomorphism.
\end{itemize}

We remark that a map $\Phi$ satisfying the two hypotheses is a local 
isometry. Now, the subgroups of the right-angled Artin group 
$\AG\cong\pi_1(Y)$ for which we shall solve the conjugacy problem are 
the fundamental groups $\pi_1(X)$ of cubical complexes~$X$ which admit
a cubical map~$\Phi\co X\to Y$ with the convexity and injectivity
property.

\begin{rem} If~$X$ and~$Y$ are both known to be $CAT(0)$ cube complexes 
then the convexity property implies the injectivity 
property -- cf.~\cite{CW1}, Theorem~1 and the remark following. 
Conversely, the two conditions, together with the knowledge that $Y$
is $CAT(0)$, imply that the complex $X$ is itself $CAT(0)$.
\end{rem}

The reader unfamiliar with the geometrical language used in stating
the conditions should remember that the convexity and injectivity 
properties are satisfied by all the subgroups of right-angled Artin 
groups discussed in Theorem~1 of~\cite{CW1}. So some typical examples
to keep in mind are those given in this paper. More generally, in order 
to get a mental image of the class of subgroups satisfying the two 
hypotheses, one can think of a subgroup whose Cayley graph sits in the 
Cayley graph of~$\AG$ in a ``flat'' way. 
Moreover, as proven by Haglund and Wise (\cite{HagWise}, Theorem 4.2), 
for a cubed complex $X$, the property of admitting map
$\Phi$ to a RAAG with the convexity and injectivity property can be
characterized purely in terms of certain combinatorial conditions
on the complex $X$ -- they call such complexes \emph{special}. 

{\bf General Notation and Conventions for the rest of the section}
\begin{itemize}
\item We fix once and for all a right-angled Artin group~$\AG$ given by 
a presentation with generators~$a_1,\ldots,a_N$, and we denote by~$Y$ 
the cubed complex associated with~$\AG$. We fix an orientation on
every edge of $Y$ and identify~$\AG$ with~$\pi_1(Y)$, using the chosen 
orientations. 
\item  We also fix a finite cubed complex~$X$ and~$\Phi\co X\to Y$ 
a cubical map satisfying the convexity and injectivity condition.
Finally, we fix a label $x_1, x_2, x_3,\ldots$ for each vertex of $X$.
\end{itemize}

Roughly speaking, our main result is the following\medskip

\emph{Claim:  Using the General Notation and Conventions of this section,
the conjugacy problem in the group $\pi_1(X)$, with respect to
any finite set of generators of $\pi_1(X)$, is solvable in linear-time.}
\medskip

Phrased in this way, however, this statement is somewhat dissatisfying, 
because we have not even stated how the generators of $\pi_1(X)$ are 
specified. A more precise statement will be given in 
Theorem~\ref{T:MainSubgrps} below.

In fact, we will not directly solve the conjugacy problem in the 
fundamental group of $X$, but a more general  
problem, namely the conjugacy problem in the fundamental \emph{groupoid} 
of $X$, in linear-time. First, we explain what precisely that means.

Let us fix a (positive) orientation for each edge of the 
complex $X$ by pulling back along $\Phi$ the orientation of edges in $Y$. 
An element of the fundamental groupoid is, 
by definition, a homotopy class of paths (with fixed endpoints) from 
some vertex $x_i$ to some vertex $x_j$. Such an element of the fundamental
groupoid can be represented by a finite sequence of successive directed 
edges, which may be traversed in the positive or in the negative direction. 
We shall call such a sequence an \emph{edge path} from $x_i$ to $x_j$.
Similarly in $Y$ we have an analogue notion of an edge path as a 
homotopy class of path specified as a sequence of positively or negatively 
directed edges.

We shall use the following very convenient way of coding edge paths 
in~$X$ and~$Y$: in $Y$, we shall simply identify closed edge paths with 
words in the letters~$a_1^{\pm 1},\ldots,a_N^{\pm 1}$.
As for~$X$, the map $\Phi$ gives rise to a coding of 
edge paths in $X$ by \emph{based words}.

\begin{defn}\label{D:basedword}
A \emph{based word} is a word of the form $x_i w x_j$, where $x_i$
and~$x_j$ are vertices of~$X$, and~$w$ is the image under $\Phi$ 
of an edge path in~$X$ starting at $x_i$ and ending at $x_j$.
The vertex $x_i$ is called the \emph{base vertex} of the based word.
\end{defn}

In other words, the edge path $x_i w x_j$ is by definition the pullback
to $X$ of the path $w$ in $Y$ which starts at $x_i$ and ends at $x_j$. 
Notice that not every word of the form $x_i w x_j$, with $x_i$ and $x_j$
vertices of $X$ and $w$ a word with letters in 
$\{a_1^{\pm 1},\ldots,a_N^{\pm 1}\}$, is a based word. However, when 
it is, then it \emph{uniquely} determines an edge path in $X$, because 
of the injectivity property. For instance, if $x_i w x_j$ is a 
based word, and if the word $w$ can be written as a concatenation
$w=w_1 w_2$, then there exists a unique vertex $x_k$ such that
$x_i w_1 x_k$ and $x_k w_2 x_j$ are based words.
For an example of based words, see again Figure~\ref{F:counterex}.

Two elements of the fundamental groupoid of $X$ can be multiplied
if the terminal vertex of the first coincides with the initial vertex
of the second. In terms of based words, $(x_i w_1 x_j)\cdot
(x_j w_2 x_k) = x_i w_1 w_2 x_k$. Two loops in $X$ are freely homotopic
if and only if they represent conjugate elements of the fundamental
groupoid. If the loops are represented by based words $x_1 w x_1$ and 
$x_2 v x_2$, then this equivalent to the existence of a based word 
$x_1 u x_2$ such that the elements of the fundamental groupoid 
represented by $x_1 u v u^{-1} x_1$ and $x_1 w x_1$ coincide.

Our main result can now be stated precisely. The proof will 
occupy the whole rest of the paper:

\begin{thm}\label{T:MainSubgrps} 
Using the General Notation and Conventions of this section, 
given two based words $x_1 w x_1$ and $x_2 v x_2$, one can decide 
whether they represent freely homotopic loops in $X$. Moreover, if 
$w$ and $v$ have length $\ell_1$ and $\ell_2$, respectively, the 
decision can be performed by an algorithm which takes time 
$O(\ell_1+\ell_2)$ on a RAM machine, where the linear constants 
depend on $X$, $Y$ and $\Phi$ only.
\end{thm}

\subsection{Base points and homotopies in the cubical complex $X$}

Why did we pass to the fundamental groupoid, rather than sticking to
the fundamental group? In other words, why do we pay so much attention
to basepoint issues? By the way of motivation, let us look at a wrong
``proof'' of Theorem~\ref{T:MainSubgrps}, and see how how
we get into trouble if we don't make basepoints explicit at every step.

{\bf Wrong Claim } Let~$\alpha$,~$\beta$ be two closed edge paths in~$X$ 
based at a common vertex~$x$. Then the loops~$\alpha$ and~$\beta$ 
represent conjugate elements of $\pi_1(X)$ if and only if 
the words~$\Phi(\alpha)$ and~$\Phi(\beta)$ represent conjugate 
elements of~$\AG$.

{\bf Wrong proof of the Wrong Claim } The implication ``$\Rightarrow$'' is 
obvious. For ``$\Leftarrow$'', we suppose that the words 
$\Phi(\alpha)$ and~$\Phi(\beta)$ represent conjugate elements of 
$\pi_1(Y)$, so the loops~$\Phi(\alpha)$ and~$\Phi(\beta)$
in~$Y$ are freely homotopic. Thus we can apply sequences of free
reductions, cyclings, and commutation relations (homotopies across
squares) in~$Y$ to each of the two loops so as to transform both of
them into some loop~$\Gamma$ in~$Y$. By the injectivity- and
convexity hypothesis, these transformations can be pulled back to
free homotopies of the original loops~$\alpha$ and~$\beta$ in~$X$.
\emph{Therefore}~$\alpha$ and~$\beta$ are both freely homotopic
to some loop~$\gamma$ in~$X$, i.e.\ they are freely homotopic. 
\hfill$\square$

This proof is almost correct, and our real proof of 
Theorem~\ref{T:MainSubgrps} shall follow this outline.
The mistake, however, is the conclusion in the very last sentence:
we can only conclude that~$\alpha$ and~$\beta$ are freely homotopic
to some loops~$\gamma$ and~$\gamma'$, respectively, where 
$\Phi(\gamma)=\Gamma=\Phi(\gamma')$. Intuitively, the loops~$\gamma$ 
and~$\gamma'$ in~$X$ may look like two different ``liftings'' 
of~$\Gamma$, we did not pay attention to basepoints!

\begin{figure}[htb]
\centerline{\input{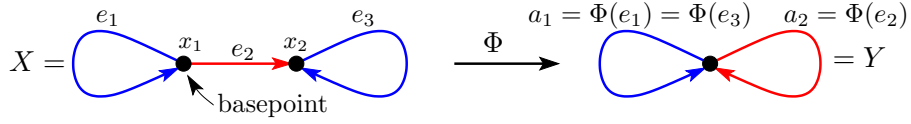}}
\caption{$\pi_1(Y)$ is the free group on two generators~$a_1=\Phi(e_1)
=\Phi(e_3)$ and~$a_2=\Phi(e_3)$. 
The loops~$e_1$ and~$e_2 e_3 e_2\inv$ are not conjugate as elements of
$\pi_1(X)$, whereas their images in~$\pi_1(Y)$ are. In order to 
describe the loops in~$X$ it is better to use the based words 
$x_1 a_1 x_1$ and $x_1 a_2 a_1 a_2\inv x_1$. The latter is conjugate
to $x_2 a_1 x_2$.}
\label{F:counterex}
\end{figure}

An explicit counterexample to the Wrong Claim illustrating the base point
problem is given in Figure~\ref{F:counterex}.

In order to prepare the proof of Theorem~\ref{T:MainSubgrps}, let us
study what homotopies of paths in~$X$ look like.

If~$\alpha$ is an 
edge path in~$X$ giving rise to a based word~$x_iwx_j$, and 
if~$x_i\wtil wx_j$ is a based word obtained from~$x_iwx_j$ by one 
application of a commutation relation (corresponding to a homotopy of a 
path in~$Y$ across a square) then there exists an edge 
path~$\wtil \alpha$ in~$X$, starting from the same vertex as~$\alpha$ 
and homotopic to~$\alpha$, which gives rise to the based 
word~$x_i\wtil wx_j$ -- this is an immediate consequence of the convexity 
condition. Similarly, free cancellations in~$w$ can be realised
by cancellations of backtracking path segments in~$\alpha$. 

Let us summarize the situation in even more geometric language.
Given a vertex~$x$ of~$X$, it is in general not true that every loop 
in~$Y$ is the image under~$\Phi$ of a path in~$X$ starting at~$x$. 
However, when such a pullback of the loop exists, then it is unique. 
Moreover, in that case 
all homotopies of the loop in~$Y$, except length-increasing ones, can 
be pulled back to based homotopies of the path in~$X$.

Let us now look more generally at \emph{free} homotopies of loops in~$X$,
i.e., homotopies that move the basepoint.

\begin{defn} Suppose that~$x w x$ is a based word. 
A \emph{parallel transport} of~$xwx$ is a replacement of the vertex~$x$ 
by a vertex~$x'$, where~$x'$ is obtained from~$x$ 
by walking along an oriented edge~$e$ with the property that the 
element~$\Phi(e)$ of~$\AG$ commutes with all the generators of~$\AG$ 
occurring in the word~$w$.
\end{defn}

Geometrically, this move corresponds to replacing a closed path based
at~$x$ by a parallel one based at~$x'$, where~$x$ and~$x'$ are joined
by an edge~$e$. The two paths together bound an annulus-shaped region 
of~$X$. Notice that, under~$\Phi$, the two paths have the same 
image~$w$ in~$Y$. 
Another way of moving the basepoint of a loop is to push it along the loop:

\begin{defn}Suppose that~$w=x y_1 y_2\ldots y_\ell x$ is a based word, 
and denote by~$e$ the unique edge of~$X$ that has one of its extremities 
equal to $x$ and such that~$\Phi(e)=y_1$. A \emph{based cycling} of the 
based word $w$ is its replacement by the word~$x' y_2\ldots y_\ell y_1 x'$, 
such that the vertex~$x'$ is the second extremity of the edge~$e$. 
\end{defn}

Geometrically, if~$\alpha$ is a loop in~$X$ based at a vertex~$x$, and 
described by a based word~$x y_1 y_2\ldots y_\ell x$, and if we apply a 
cycling operation (in the sense of section~\ref{S:RAAGAlg}) to the 
word~$y_1 y_2\ldots y_\ell$, 
then this cycling can be pulled back to~$X$ to a based cycling of the
based word, yielding a loop~$\widetilde{\alpha}$, which looks exactly 
like~$\alpha$, except that it based at a different vertex~$x'$, ``one 
notch further along the loop''.

\begin{exam}\label{E:counterex} 
In the example of Figure \ref{F:counterex}, we can apply a based 
cycling to the based word~$x_1 a_2 a_1 a_2\inv x_1$, yielding
$x_2 a_1 a_2\inv a_2 x_2$. After a cancellation, we obtain the based
word~$\gamma_2=x_2 a_1 x_2$. We note that this is different from the based 
word~$\gamma_1=x_1 a_1 x_1$,
which was also discussed in that example -- in fact, the based words
$x_2 a_1 x_2$ and~$x_1 a_1 x_1$ are not even related by parallel transport
(because~$\Phi(e_2)=a_2$ does not commute with~$a_1$). As we shall see
in Lemma~\ref{L:SameImage}, this implies that  
the two loops~$e_1$ and~$e_2 e_3 e_2\inv$ are not freely 
homotopic in~$X$. 
\end{exam}

Also note that a cyclic reduction of a word on the generators of~$A$ 
and their inverses can be decomposed as a cycling, followed by a usual 
cancellation of letters, and each of these operations can be pulled back 
to operations on the loop in~$X$.
Summarizing the last few paragraphs, we have the following

\begin{keyobs}\label{KO:pullback}
If~$\alpha$ is a loop in~$X$
then all non-length-increasing free 
homotopies of the loop~$\Phi(\alpha)$ in~$Y$ 
can be pulled back to free homotopies of~$\alpha$. 
Thus for a based word $x w x$, all cancellations, applications of 
commutation relations, cyclings, and cyclic reductions of the word~$w$
can be pulled back to analogue cancellations, commutation relations, 
and based cyclings of the based word. 
Similarly, if $x_1 w x_2$ is a based word,
and if the word $w$ can be transformed into a word $w'$ by applying
cancellations and commutation relations, then $x_1 w' x_2$ is again
a based word.
\end{keyobs}

\subsection{The linear-time solution to the conjugacy problem}

The aim of this subsection is to prove Theorem~\ref{T:MainSubgrps}. 
We recall that we are considering two based words $x_1 w x_1$ and
$x_2 v x_2$ representing two loops in $X$ traversing $\ell_1$ and
$\ell_2$ edges, respectively. A necessary condition for these loops
being conjugate in the fundamental groupoid of $X$ is that the 
words~$w_1$ and~$w_2$ represent conjugate elements of the right-angled 
Artin group $\AG$. In geometric terms, for the two loops to be freely 
homotopic in $X$, their images under~$\Phi$ in~$Y$ must be freely homotopic.
This is a condition which we can check in time $O(\ell_1+\ell_2)$
by the results of Section~\ref{S:RAAGAlg}. However, this condition is
not sufficient, as seen in Example~\ref{E:counterex}. So let us now try 
to refine this approach.

\begin{prop}\label{L:ReduceToSameProj}
There is an algorithm with running time $O(\ell_1+\ell_2)$ whose input
consists of two based words $x_1 w x_1$ and 
$x_2 v x_2$ of lengths~$\ell_1$ and~$\ell_2$, and which outputs
\begin{enumerate}
\item either the information that they do not represent freely homotopic
loops in~$X$, or
\item two based words $x'_1 \wtil w_1 \ldots \wtil w_k x'_1$ and 
$x_2' \wtil w_1 \ldots \wtil w_k x_2'$, representing two loops in~$X$ 
which are respectively freely homotopic to the original two, and where
the $\wtil w_i$ are mutually commuting cyclic normal forms.
\end{enumerate}
\end{prop}

\begin{proof}[Proof of Proposition~\ref{L:ReduceToSameProj}]
As seen in Section~\ref{S:RAAGAlg} we can decide in linear-time 
whether $w$ and $v$ represent conjugate
elements of $\AG$. If they do not, then the two based words do not
represent conjugate elements of the fundamental groupoid either, 
and it suffices to output this information (case (1)). 

For the rest of the proof we have to deal with the case where $w$ and $v$
do represent conjugate elements of~$\AG$.

We already know from Section~\ref{S:RAAGAlg} that the word $w$ can, 
by a sequence of cancellations, commutation relations and cyclings be 
transformed into a word $w'$  with the required decomposition
$w'=w'_1\ldots w'_k$.
Moreover, we know how to calculate the word $w'$ in linear-time.

We also know from the Key Observation~\ref{KO:pullback} above that the 
transformation of the word $w$ into the word~$w'$ 
can be pulled back to a transformation of the 
based word~$x_1 w x_1$ into a based word~$x_3 w' x_3$. 
Our next task is to determine the corresponing base 
vertex~$x_3$ in linear-time. 

We shall fulfill this task by ``carrying along information about the 
base vertex in~$X$ during the algorithm''.
While running the algorithm of Section~\ref{S:RAAGAlg}, the only steps
that affect the base vertex are the cyclings of pilings (including cyclic
reductions of pilings, which can be decomposed as cyclings, followed by 
cancellations of tiles): when we cycle an~$a_j^\pm$-tile, we have to 
determine how the base vertex is affected. However, this can be done
simply by a lookup in a finite, precalculated list: for every vertex~$x$ 
of~$X$, for every generator~$a_j$ of~$A$, and for every 
$\epsilon\in\{-1,1\}$, this list must tell us at which vertex of~$X$ 
we arrive if we pull back the loop~$a_j^\epsilon\in \AG=\pi_1(Y)$
to a path in~$X$ starting at~$x$ (if that is possible). 
Since the algorithm of Section~\ref{S:RAAGAlg} performs a linearly 
bounded number of cyclings, we can calculate the new base vertex~$x_3$ 
in time~$O(\ell_1)$.

In a similar manner we can algorithmically transform the based word 
$x_2 v x_2$ into a word $x'_2 \wtil w x'_2$, where $\wtil w$ is equipped 
with an analogue decomposition 
${\nobreak \wtil w=\wtil w_1\cdots \wtil w_j}$.

But since $w$ and $v$ represented conjugate elements of~$\AG$,
we have, by the results of Section~\ref{S:RAAGAlg}, that the words~$w'$
and~$\wtil w$ are in fact the same, at least after a reordering of the factors of~$w'$ and a linearly bounded number of cyclings of each 
factor $w'_i$; in particular, we have $j=k$. 
Moreover, the Boyer-Moore algorithm tells us how many letters from each 
factor we have to cycle in order to achieve this. Thus we can transform the
based word $x_3 w' x_3$ into the based word $x'_1 \wtil w x'_1$ for some 
vertex $x'_1$, using a reordering of the factors (which does not 
affect the base vertex) and a linearly bounded number of based cyclings.
\end{proof}

Thus in order to prove Theorem~\ref{T:MainSubgrps}, it is enough to 
prove it for the special case $v=w=\wtil w_1\ldots \wtil w_k$, where the 
words $\wtil w_1,\ldots,\wtil w_k$ are mutually commuting cyclic normal 
forms. (For instance, this is the situation of Example~\ref{E:counterex}, 
where we need to decide if the based words $x_1 a_1 x_1$
and $x_2 a_1 x_2$ represent freely homotopic loops in~$X$.)
For the rest of the proof of Theorem~\ref{T:MainSubgrps} we fix 
such a word~$\wtil w$, with such a decomposition.

Suppose a based word $x_1 \wtil{u} x_2$ 
is such that $x_1 \wtil{u} \wtil w \wtil{u}^{-1} x_1$ and $x_1 \wtil w x_1$ 
represent the same element of the fundamental groupoid.
Then in particular the elements of~$\AG$ represented 
by $\wtil{u}$ and $\wtil w$ commute: we have 
$[\wtil u][\wtil w][\wtil u]^{-1}=[\wtil w]$ in $\AG$. 

As seen in Section~\ref{SS:centralizer}, and using the notation of 
this section, the word $\wtil{u}$ is equivalent to another word $u$ of 
the form
$$
u=z_1^{p_1}\ldots z_k^{p_k}\,\zeta
$$
where $p_1,\ldots,p_k$ are integers and $\zeta$ is a word whose letters
are generators of~$\AG$ which commute with, but are different from, 
all the generators occurring in~$w$, and their inverses. We shall call 
such a word $u$ a word in 
\emph{preferred form}. We define the \emph{norm} $\| u\|$ of~$u$ by
$$
\|u\|=\sum_{i=1}^k |p_i| + \mathrm{length}(\zeta)
$$

We are now ready state an algorithmically checkable criterion for 
$x_1 \wtil w x_1$ and $x_2 \wtil w x_2$ representing conjugate elements 
(i.e.\ representing freely homotopic loops in~$X$):

\begin{lemma}\label{L:SameImage}
The two based words $x_1 \wtil w x_1$ and $x_2 \wtil w x_2$ represent 
conjugate elements in the fundamental groupoid if and only 
if there exists a based word~$x_1 u x_2$
such that~$u$ is a word in preferred form with
\begin{equation}\label{F:ExpSum}
\|u\|\leqslant \#\{\textrm{vertices\ of }X\}
\end{equation}
\end{lemma}

\begin{proof}
We first suppose that an edge path $x_1 u x_2$ exists, where $u$ is a 
word in preferred form.
Then the word $u \wtil w u^{-1}$ can be transformed into the word 
$\wtil w$ by a
finite number of commutation relations and cancellations (but no
length-increasing transformations). By Key Observation~\ref{KO:pullback}, 
this homotopy can be pulled back to $X$, to yield a based homotopy 
between the paths in~$X$ represented by the based words 
$x_1 u \wtil w u^{-1} x_1$ and $x_1 \wtil w x_1$.
In other words, the elements $x_1 \wtil w x_1$ and $x_2 \wtil w x_2$ 
are conjugate, with conjugating element $x_1 u x_2$.

Conversely, let us suppose that a conjugating element in the fundamental 
groupoid exists, and is represented by a based word $x_1 \wtil{u} x_2$. 
This means that there exists an edge path in $X$ from $x_1$ to~$x_2$ 
such that reading out the edge labels along the path yields the 
word~$\wtil{u}$. 
As seen before, $[\wtil{u}]$ belongs to 
the subgroup of~$\AG$ generated the elements $[z_1],\ldots ,[z_k]$ and 
$[a_{j_1}],\ldots,[a_{j_m}]$. Thus there is a word $u$ in preferred 
form which can be obtained from $\wtil{u}$ by a sequence of reductions and
commutation relations. By Key Observation~\ref{KO:pullback}, 
$x_1 u x_2$ is also a based word, i.e.\ it also represents an edge path 
in~$X$.

We have shown the existence of a based word $x_1 u x_2$ with~$u$
a word in preferred form, and without loss of generality we can
suppose that $u$ is chosen so that $\|u\|$ is minimal among all such
based words.

Now for $t$ in $\{ 0,\ldots,\|u\|\}$ let us denote by $x(t)$ the vertex 
of~$X$ obtained by a walk in $X$ starting at $x_1$ and following the edges 
of~$X$ according to the~$t$ first subwords. Now, if this function
$$
\{0,1,\ldots,\|u\|\} \longrightarrow \{\textrm{vertices 
of }X\}\ , \ \ \ t \mapsto x(t)
$$
is not injective (for instance, if~$\|u\|$ is larger than the number
of vertices of~$X$), then there exists a strictly shorter edge path in
$X$ represented by a based word $x_1 u' x_2$ with $u'$ also in
preferred form, obtained by cutting out some segment of the previous edge 
path (c.f.\ the paragraph following Definition~\ref{D:basedword}). 
This is in contradiction to the choice of $u$, and we can conclude
that we have $\|u\|\leqslant \#\{\textrm{vertices\ of }X\}$
\end{proof}

Let us now prove that the condition of Lemma~\ref{L:SameImage} can be 
checked algorithmically in linear-time, i.e.~in time 
$O(\ell)$, where~$\ell$ is the length of the word $w$. 

Firstly, recalling that the centralizer of $[\wtil w]$ is generated by
a finite number of elements (some of them represented by the words $z_1,\ldots,z_k$ and the others equal to certain generators of~$\AG$),
we observe that there is a universal upper bound on the number of
generators, namely the number of generators of~$\AG$. Moreover,
as seen in Proposition~\ref{P:commutator}, words representing these 
generators can be determined in linear time. 

Now there is a very simpleminded linear-time algorithm to check for the
existence of a conjugating element: for \emph{all} words $u$ 
in preferred form satisfying condition~(\ref{F:ExpSum}) check whether 
$x_1 u x_2$ is a based word, i.e.\ whether there exists an edge path 
in~$X$ represented by the based word $x_1 u x_2$. Indeed, there is a 
universal bound on the number of words to be checked, and for each 
word~$u$ the check takes linear time (since the length of the 
words~$z_i$ can grow linearly with the length of $\wtil w$).

Here is a summary of the whole algorithm: given two based words~$x_* w x_*$ 
and $x_* v x_*$ representing loops $\alpha$ and $\beta$ in $X$, 

(1) Apply steps (i) and (ii) of the algorithm of 
Section~\ref{SS:CyclicNF}, always carrying 
along the base vertex, to find graphs~$\Delta_j(w)$ ($j=1,\ldots,k$),  
$\Delta_j(v)$ ($j=1,\ldots,k'$), base vertices~$x_1$,~$x_2$, and 
based words $x_1w_1\ldots w_kx_1$ 
and $x_2v_1\ldots v_{k'}x_2$ representing loops that are
freely homotopic to~$\alpha$ and~$\beta$.

(2) If~$k\neq k'$, or if the collections of full subgraphs~$\Delta_j(w)$ 
and $\Delta_j(v) \subset \graph$ are not the same, or if for 
some~$j$ between $1$ and $k$ the words $v_j$ and $w_j$ do not have the 
same length~$\ell_j$, return ``NO''.

(3) Apply step (iii)(a) of the algorithm of 
Section~\ref{SS:CyclicNF} to each of the $k$ factors, 
always carrying along the base vertices, to transform
$x_1w_1\ldots w_kx_1$ into a based word $x_3 w'_1\ldots w'_k x_3$ 
and similarly $x_2 v_1\ldots v_k x_2$ into 
$x'_2\wtil w_1\ldots \wtil w_{k}x'_2$,
where all words $w'_i$ and $\wtil w_i$ are cyclic normal forms.

(4) For each factor, use a standard pattern matching algorithm to decide 
if $w'_i=\wtil w_i$ as cyclic words. If no, return ``NO''. If yes, keep 
in mind how many cyclings of each factor $w'_i$ are required to achieve 
equality $w'_i=\wtil w_i$ as (non-cyclic) words.

(5) Perform the required \emph{based} cyclings of 
$x_3w'_1\ldots w'_k x_3$ to obtain a based word of the form
$x'_1\wtil w_1\ldots \wtil w_k x'_1$.

(6) Calculate the minimal roots $z_i$ of the words $\wtil w_i$, as
explained in Section~\ref{SS:centralizer}. Also determine
the set of generators that commute with all the letters occurring in
the words $\wtil w_i$, but do not occur in any of them.

(7) Check, for all words $u$ in preferred from satisfying 
condition~(\ref{F:ExpSum}), whether there exists an edge path in $X$
represented by the based word $x'_1 u x'_2$. If for one of the 
words~$u$ the answer is affirmative, then return ``YES''.
Otherwise return ``NO''.

{\bf Acknowledgement } We thank Sam Sang-Hyun Kim, Tim Hsu, 
Lucas Sabalka, and Michah Sageev for interesting conversations.

%|<------------------------------------------------------------------------>|

\end{document}